\documentclass[11pt]{article}
\usepackage[T1]{fontenc}
\usepackage{lmodern,amsmath,amssymb,amsthm,mathtools,booktabs,microtype}
\usepackage[a4paper,hmargin=28mm,vmargin=25mm]{geometry}
\usepackage{xcolor}
\usepackage{algorithm,algpseudocode}
\usepackage[colorlinks=true,linkcolor=blue!50!black,citecolor=blue!50!black,urlcolor=blue!50!black]{hyperref}
\newtheorem{theorem}{Theorem}[section]
\newtheorem{proposition}[theorem]{Proposition}
\newtheorem{lemma}[theorem]{Lemma}
\theoremstyle{remark}\newtheorem{remark}[theorem]{Remark}
\newcommand{\Sn}{\mathbb S^n} \newcommand{\cK}{\mathcal K}
\newcommand{\ip}[2]{\langle #1,#2\rangle}
\newcommand{\fn}[1]{\lVert #1\rVert_F} \newcommand{\E}{\mathbb E}
\DeclareMathOperator{\ran}{ran} \DeclareMathOperator{\tr}{tr}
\newcommand{\qf}{q}
\newcommand{\loss}{\mathcal E}
\hypersetup{pdftitle={A note on approximate projection onto the semidefinite complementarity set},pdfauthor={Chao Ding}}
\begin{document}
	\begin{center}
		{\LARGE A note on approximate projection onto the\\[2pt]
			semidefinite complementarity set\par}
		\medskip
		{\large Chao Ding\par}
		Academy of Mathematics and Systems Science, Chinese Academy of Sciences\par
		Beijing, China\par
		\href{mailto:dingchao@amss.ac.cn}{\nolinkurl{dingchao@amss.ac.cn}}\par
		\smallskip
		25 November 2026\par
	\end{center}
	\smallskip
	\begin{abstract}
		The proximal slack update in an augmented Lagrangian method for semidefinite complementarity constraints requires a joint matrix projection. Given real symmetric matrices $A$ and $B$, we seek the nearest pair $(X,Z)$ with $X$ positive semidefinite, $Z$ negative semidefinite, and $XZ=0$. Hou, Peng and Ding obtain a squared-loss factor of $2$ by fitting the pair in an eigenbasis of $A+B$. To improve this bound, we replace the orthogonal projector onto the subspace allowed for $X$ by a symmetric matrix $T$ with eigenvalues in $[0,1]$. From an optimal solution of the resulting semidefinite program, we form candidate subspaces by starting at the zero subspace and adding the eigenspaces of $T$ in decreasing order of eigenvalue. For each candidate, we compress $A$ to that subspace and $B$ to its orthogonal complement. We keep the positive eigenvalues of the first compression and the negative eigenvalues of the second. The best pair has squared loss at most $1.189$ times the minimum for arbitrary inputs and at most $1.15487$ times the minimum when $A$ is positive semidefinite and $B$ is negative semidefinite. For arbitrary inputs the factor improves to $1.18$ if $T$ has at most two distinct eigenvalues in $(0,1)$. The threshold-rounding proofs average over a uniform spectral threshold and use computer-assisted bounds for the resulting scalar inequalities.
	\end{abstract}
	
	\section{The projection problem}
	Let $\Sn$ denote the real symmetric $n\times n$ matrices and $I$ the identity matrix. For real matrices of the same size we use the Frobenius inner product $\ip{H}{K}=\tr(H^\top K)$ and norm $\fn{H}=\sqrt{\ip{H}{H}}$. For symmetric matrices, $H\preceq K$ means that $K-H$ is positive semidefinite. For $H\succeq0$, $H^{1/2}$ denotes its positive semidefinite square root. We write $\ran H$ and $\ker H$ for the range and null space of $H$. Given $A,B\in\Sn$, we seek the squared distance from $(A,B)$ to the graph of the normal-cone mapping of the positive semidefinite cone,
	\[
	\{(X,Z)\in\Sn\times\Sn:X\succeq0,\ Z\preceq0,\ XZ=0\}.
	\]
	A matrix $Z$ is normal to the positive semidefinite cone at $X\succeq0$ if and only if $\ip{Z}{H-X}\le0$ for every $H\succeq0$. Taking $H=0$ and $H=2X$ gives $\ip{Z}{X}=0$. Taking $H=vv^\top$ then gives $v^\top Zv\le0$ for every $v\in\mathbb R^n$, so $Z\preceq0$. Since $-\ip{Z}{X}=\fn{X^{1/2}(-Z)^{1/2}}^2$, orthogonality of $X$ and $Z$ is equivalent to $XZ=0$ under these sign conditions. Conversely, $Z\preceq0$ and $XZ=0$ give $\ip{Z}{H-X}=\ip{Z}{H}\le0$ for every $H\succeq0$.

	Set $Y=-Z$ and write $C=A$, $D=-B$. Reflection preserves the Frobenius distance. For a fitted pair $(X,Y)$, let $\loss(X,Y)$ denote its squared loss and $v_*$ the minimum over complementary pairs. Thus
	\begin{equation}\label{eq:problem}
		\loss(X,Y)=\fn{C-X}^2+\fn{D-Y}^2,\qquad
		v_*=\min_{X,Y\succeq0,\ XY=0}\loss(X,Y).
	\end{equation}
	A feasible pair $(X,Y)$ in \eqref{eq:problem} gives $(X,-Y)$ in the original variables. The case $C,D\succeq0$ corresponds to $A\succeq0$ and $B\preceq0$. The ranges of $C$ and $D$ may overlap, and the matrices need not commute. For $H\in\Sn$, let $H_+$ denote its positive part, obtained by replacing each eigenvalue $\lambda$ by $\max\{\lambda,0\}$ in an orthonormal eigenbasis of $H$.
	
	For a primal--dual feasible pair in a linear semidefinite program, the duality gap is the trace inner product of the primal matrix and the dual slack. Both matrices are positive semidefinite, so a zero gap requires their ranges to be orthogonal. Wu et al.~\cite{WuZhangZhangLu2016} use this complementarity condition in an inverse problem that adjusts the objective and right-hand-side data of a linear semidefinite program. A rank bound gives another example. Ding, Sun and Ye~\cite[eqs.~(2)--(3)]{DingSunYe2014} use a complementary matrix to impose a rank bound on a correlation matrix $X$. For an integer $0\le r\le n$, the condition $\operatorname{rank}X\le r$ is equivalent to the existence of $Y$ with $0\preceq Y\preceq I$, $\tr Y=n-r$, and $XY=0$. Indeed, if $\operatorname{rank}X\le r$, we can take $Y$ to project onto an $(n-r)$-dimensional subspace of $\ker X$. Conversely, any such $Y$ has $\operatorname{rank}Y\ge\tr Y=n-r$ and $\ran Y\subseteq\ker X$.

	The projection itself occurs in the augmented Lagrangian method of Hou, Peng and Ding~\cite[eqs.~(33)--(34), (57)]{HouPengDing2026}. The method introduces a separate pair of slack matrices constrained to be complementary. With the other variables fixed, the slack objective contains a squared distance to the multiplier-shifted constraint values and a proximal squared distance to the previous slack pair. If their weights are $\rho/2$ and $b/2$, respectively, with $\rho,b>0$, completing the square gives $(\rho+b)\loss(X,Y)/2$ up to an additive constant. Here $(C,D)$ is the weighted average of the two centers in the reflected variables. Thus an exact slack update solves \eqref{eq:problem}. A feasible pair of loss at most $c v_*$, for $c\ge1$, has slack objective at most $(\rho+b)(c-1)v_*/2$ above the exact block minimum. Such a pair can be used as a trial update in the inexact method, which also tests descent and a limiting-normal residual before acceptance~\cite[Algorithm~2, eqs.~(39)--(40)]{HouPengDing2026}.

	When $n$ is part of the input, deciding whether $v_*$ is at most a given rational threshold is NP-hard, even for rational positive definite $C$ and $D$ \cite[Remark~4 and Appendix~E]{HouPengDing2026}.

	\begin{remark}[No fully polynomial-time approximation scheme]\label{rem:no-fptas}
		A fully polynomial-time approximation scheme (FPTAS) for \eqref{eq:problem} would, for any rational $0<\varepsilon<1$, return a feasible pair of loss $\widehat v\le(1+\varepsilon)v_*$ in time polynomial in the input bit length and $1/\varepsilon$. Even for rational $C,D\succ0$, such a scheme would imply $\mathrm P=\mathrm{NP}$. To see this, take a simple unweighted graph on $n\ge2$ vertices and let $m$ be its maximum cut size. The reduction in \cite[Appendix~E]{HouPengDing2026} constructs the matrices in polynomial time. With $\delta=(1024n^{16})^{-1}$ and $c_0=\fn{D}^2$, its final estimate is $\big|v_*-c_0+8\delta^2m\big|\le\tfrac12\delta^2$. The diagonal entries of $D$ are $1^2,\ldots,n^2$, and its off-diagonal part has Frobenius norm less than one. Hence $c_0\le n^5+1$. The feasible pair $(C,0)$ gives $v_*\le c_0$. Choose $\varepsilon=\delta^2/c_0$, so that $1/\varepsilon=O(n^{37})$ and $0\le\widehat v-v_*\le\delta^2$. It follows that $(c_0-\widehat v)/(8\delta^2)$ differs from $m$ by at most $3/16$. Rounding to the nearest integer would recover the maximum cut size in polynomial time and hence solve the NP-complete Simple Max-Cut decision problem~\cite{GareyJohnsonStockmeyer1976}. The same conclusion holds for the Frobenius distance itself, since a distance factor of $1+\varepsilon/3$ gives a squared-loss factor at most $1+\varepsilon$ for $0<\varepsilon<1$.
	\end{remark}

	Let $S=C-D=A+B$. Its spectral split $(S_+,(-S)_+)$ is feasible for \eqref{eq:problem}. In any orthonormal eigenbasis of $S$, let $V_\Sigma$ be the least loss among complementary pairs diagonal in that basis. On each basis vector $e$, the diagonal fit keeps the larger of $(e^\top C e)_+$ and $(e^\top D e)_+$ in the corresponding matrix and sets the other entry to zero, choosing the first matrix in a tie. Hou, Peng and Ding~\cite[Proposition~8, eq.~(77), and Appendix~B]{HouPengDing2026Revised} give
	\[
		V_\Sigma \le \loss(S_+,(-S)_+) = 2\fn{C-S_+}^2 \le 2v_* .
	\]
	The spectral split attains the factor $2$ when $C=D=[1]$, for which its loss is $2$ and $v_*=1$. Taking the best diagonal fit over chosen eigenbases of $C$, $D$, and $C-D$ also gives loss at most $2v_*$, since the third basis alone guarantees this bound.

	The constraint $XY=0$ requires the ranges of $X$ and $Y$ to be orthogonal. Projecting $C$ and $D$ separately onto the positive semidefinite cone can leave their fitted ranges overlapping. Once a subspace for $X$ and its orthogonal complement for $Y$ have been chosen, two positive semidefinite projections give the best fit. We relax the projector onto the subspace for $X$ to a matrix $0\preceq T\preceq I$. Thresholding its eigenvalues gives the candidate subspaces. Table~\ref{tab:bounds} gives the loss bounds. Whether the bounds for threshold rounding are sharp remains open (Remark~\ref{rem:tight}).
	\begin{table}[ht]
		\centering
		\caption{Squared-loss factors for \eqref{eq:problem}. The spectral split and diagonal refit use an eigenbasis of $C-D$. Threshold rounding uses an optimal solution of \eqref{eq:relax}. For $v_*>0$, the square root of each factor bounds the distance ratio.}
		\label{tab:bounds}
		\small\begin{tabular}{@{}p{.36\linewidth}p{.39\linewidth}r@{}}\toprule
			Method & Hypothesis & Squared-loss factor\\ \midrule
			Spectral split or diagonal refit & $C,D\in\Sn$ & $2$\\
			Threshold rounding & $C,D\in\Sn$ & $1.189$\\
			Threshold rounding & $C,D\succeq0$ & $1.15487$\\
			Threshold rounding & At most two distinct eigenvalues of $T$ in $(0,1)$ & $1.18$\\\bottomrule
		\end{tabular}
	\end{table}
	
	\section{Semidefinite relaxation and spectral threshold rounding}
	Let $P$ be an orthogonal projector and set $Q=I-P$. Write $v(P)$ for the minimum of $\loss(X,Y)$ subject to $X,Y\succeq0$, $X=PXP$, and $Y=QYQ$. The last two conditions place the ranges in orthogonal subspaces, so every such pair is complementary. Decomposing the inputs into blocks relative to $\ran P$ and $\ran Q$ gives
	\begin{equation}\label{eq:hard}
		(X_P,Y_P)=((PCP)_+,(QDQ)_+),\quad
		v(P)=\fn C^2+\fn D^2-\fn{(PCP)_+}^2-\fn{(QDQ)_+}^2 .
	\end{equation}
	For the first matrix, only the $P$-block can vary. By the Frobenius projection formula for the positive semidefinite cone~\cite[Theorem~2.1]{Higham1988}, truncating the negative eigenvalues of $PCP$ gives the least loss on that block. The same argument applies to $QDQ$ on $\ran Q$. These two compressed projections also appear in the fixed-subspace slack fit in \cite[Sec.~6.1]{HouPengDing2026}. Conversely, for any complementary pair one may take $P$ to be the projector onto $\ran X$, since $\ran Y\subseteq\ker X$. It follows that $v_*=\min_{P=P^\top=P^2}v(P)$. The minimum is attained because $v(P)$ is continuous on the compact set of orthogonal projectors.
	
	\subsection{The relaxation}
	At a projector $P$, the map $Z\mapsto PZP$ keeps precisely the block on $\ran P$. To extend this map to $0\preceq T\preceq I$, work on the full space $\mathbb R^{n\times n}$ with the Frobenius inner product. Define three linear operators by their action on a matrix $Z$,
	\[
	S_T(Z)=\tfrac12(TZ+ZT),\qquad J_T(Z)=TZ-ZT,\qquad H_T(Z)=2Z-TZ-ZT.
	\]
	Let $\dagger$ denote the Moore--Penrose inverse, and set $\cK_T=S_T-\tfrac12J_T H_T^\dagger J_T$. Choose an orthonormal eigenbasis $(e_i)$ of $T$. If $e_i,e_j$ have eigenvalues $x,y$, the matrix unit $e_i e_j^\top$ is a common eigenvector of $S_T,J_T,H_T$, with eigenvalues $(x+y)/2$, $x-y$, and $2-x-y$, respectively. These matrix units form an orthonormal basis of $\mathbb R^{n\times n}$, so the three operators are self-adjoint and commute. The eigenvalue of $\cK_T$ on $e_i e_j^\top$ is
	\begin{equation}\label{eq:k}
		k(x,y)=\frac{x(1-x)+y(1-y)}{2-x-y}\quad ((x,y)\ne(1,1)),\qquad k(1,1)=1.
	\end{equation}
	Since $1-k(x,y)=((1-x)^2+(1-y)^2)/(2-x-y)$ away from $(1,1)$, and $k(x,y)$ is a weighted average of $x$ and $y$ with weights $1-x,1-y$, we have $0\le k\le1$ on $[0,1]^2$. The weighted-average formula also gives continuity at $(1,1)$. Hence $\cK_T$ is positive semidefinite and all its eigenvalues lie in $[0,1]$. At a projector $P$, the multiplier is one when both eigenvalues are one and zero otherwise, so $\cK_P(Z)=PZP$. For $U\in\mathbb R^{n\times n}$, define the extended-real quadratic form
	\[
	\qf_T(U)=
	\begin{cases}
		\ip{U}{\cK_T^\dagger U},&U\in\ran\cK_T,\\
		+\infty,&\text{otherwise}.
	\end{cases}
	\]
	For a symmetric input $H$, let $\Phi_H(T,U)$ be the contribution of $U$ to the relaxed loss. We write $L(T,U,W)$ for the total and $\ell_*$ for its minimum,
	\begin{align}
		\Phi_H(T,U)&=\fn H^2-2\ip H U+\qf_T(U),\notag\\
		L(T,U,W)&=\Phi_C(T,U)+\Phi_D(I-T,W),\notag\\
		\ell_*&=\min_{0\preceq T\preceq I,\ U,W\succeq0} L(T,U,W).\label{eq:relax}
	\end{align}
	A feasible triple $(T,U,W)$ is \emph{finite} if $L(T,U,W)<\infty$, equivalently if both $\qf_T(U)$ and $\qf_{I-T}(W)$ are finite. Write $x_i$ for the eigenvalues of $T$ and use the corresponding eigenbasis for the entries $U_{ij}$. The range condition in the definition of $\qf_T$ requires $U_{ij}=0$ whenever $k(x_i,x_j)=0$. On this domain,
	\[
	\qf_T(U)=\sum_{k(x_i,x_j)>0}U_{ij}^2/k(x_i,x_j).
	\]
	For the linear matrix inequality, use the fixed standard matrix-unit basis of $\mathbb R^{n\times n}$, independently of $T$. Let $z\in\mathbb R^{n^2}$ be the vector formed by stacking the columns of $U$, using the same order for the basis. For $\tau\in\mathbb R$, the epigraph condition $\qf_T(U)\le\tau$ is equivalent to
	\begin{equation}\label{eq:lmi}
		\begin{pmatrix}\tau&z^\top&0\\z&S_T&J_T\\0&J_T&2H_T\end{pmatrix}\succeq0 .
	\end{equation}
	Simultaneous diagonalization gives $\ker H_T\subseteq\ker J_T$. By the generalized Schur complement criterion~\cite[Appendix~A.5.5]{BoydVandenberghe2004}, eliminating $2H_T$ leaves a block matrix with diagonal blocks $\tau,\cK_T$ and off-diagonal block $z$. Positive semidefiniteness requires $z$ to be orthogonal to $\ker\cK_T$. Completing the square on $\ran\cK_T$ then gives $\qf_T(U)\le\tau$. All entries in \eqref{eq:lmi} are affine in $(T,U,\tau)$. Use a second scalar variable to bound $\qf_{I-T}(W)$, with $(I-T,W)$ in place of $(T,U)$ in \eqref{eq:lmi}. Minimizing the sum of the two scalar bounds and the affine terms in \eqref{eq:relax} gives a semidefinite program with two linear matrix inequalities of order $2n^2+1$.
	When $T=P$ is a projector, $\cK_P$ projects onto the matrices $U=PUP$. On this subspace $\qf_P(U)=\fn U^2$. Similarly, $\qf_{I-P}(W)$ is finite precisely when $W=(I-P)W(I-P)$ and then equals $\fn W^2$. Minimizing the relaxed objective at $T=P$ therefore gives $v(P)$, so $\ell_*\le v_*$. For every admissible $T$, the bound $\qf_T(U)\ge\fn U^2$ gives $L(T,U,W)\ge\fn{C-U}^2+\fn{D-W}^2$. Thus $0\le\ell_*$ and $U,W$ are bounded on every finite lower level set. The epigraph in \eqref{eq:lmi} is closed and the admissible set for $T$ is compact, so the relaxation attains its minimum.
	
	\subsection{Threshold rounding}
	Let $\lambda_1>\cdots>\lambda_m$ be the distinct eigenvalues of $T$, with eigenspaces $E_1,\ldots,E_m$. Start with the zero subspace and successively add $E_1,\ldots,E_m$, ending with $\mathbb R^n$. After step $j$ the subspace contains exactly the eigenspaces with eigenvalues at least $\lambda_j$. Each repeated eigenvalue enters with its entire eigenspace. Let $\mathcal P(T)$ consist of the orthogonal projectors onto these $m+1\le n+1$ subspaces. The least loss over this family is
	\begin{equation}\label{eq:V}
		V(T)=\min_{P\in\mathcal P(T)}v(P).
	\end{equation}
	
	\begin{algorithm}[ht]
		\caption{Semidefinite relaxation and spectral threshold rounding}\label{alg:scan}
		\begin{algorithmic}[1]
			\Require Real symmetric inputs $C,D$.
			\State Choose an optimal triple $(T,U,W)$ for \eqref{eq:relax}.
			\State Order the distinct eigenvalues of $T$ and form the nested projectors $\mathcal P(T)$.
			\For{$P\in\mathcal P(T)$}
			\State Set $Q=I-P$, compute $(PCP)_+$ and $(QDQ)_+$, and evaluate $v(P)$ by \eqref{eq:hard}.
			\EndFor
			\State \Return a pair $(X_P,Y_P)$ with the least value of $v(P)$.
		\end{algorithmic}
	\end{algorithm}
	For fixed $T$, Algorithm~\ref{alg:scan} returns a pair of loss $V(T)$, even if several projectors attain the minimum. For any threshold $s\in[0,1]$, write $P_s(T)$ for the projector onto the sum of the eigenspaces of $T$ with eigenvalues at least $s$. For fixed $T$ abbreviate it to $P_s$. In the proofs $s$ is uniform on $[0,1]$, and $\E$ denotes expectation with respect to the random threshold. Every $P_s$ belongs to $\mathcal P(T)$, so $V(T)\le v(P_s)$ at each threshold and hence $V(T)\le\E\,v(P_s)$.
	
	\begin{theorem}\label{thm:main}
		For every finite feasible triple $(T,U,W)$ of \eqref{eq:relax},
		\[
		V(T)\le 1.189\, L(T,U,W)\quad(C,D\in\Sn), \qquad
		V(T)\le 1.15487\, L(T,U,W)\quad(C,D\succeq0).
		\]
		Let $c$ be the constant in the applicable inequality. Algorithm~\ref{alg:scan} returns a pair of squared loss $V(T)\le c v_*$. More generally, for an absolute tolerance $\eta\ge0$, if a finite feasible triple satisfies $L(T,U,W)\le\ell_*+\eta$, threshold rounding gives $V(T)\le c(v_*+\eta)$.
	\end{theorem}
	For an optimal triple, $L(T,U,W)=\ell_*\le v_*$. For the last claim, use $L(T,U,W)\le\ell_*+\eta\le v_*+\eta$. Here $\eta$ measures absolute error in the relaxation objective. Both claims assume exact eigenspaces and exact fits on the selected subspaces. A bit-complexity bound for relative error would also require an estimate of the cost of solving the relaxation to that accuracy. If $v_*=0$, an optimal relaxed solution gives $V(T)=0$, so the returned pair equals the input pair in \eqref{eq:problem}.
	
	\begin{remark}\label{rem:tight}
		The bounds $1.189$ and $1.15487$ in Theorem~\ref{thm:main} and $59/50$ in Proposition~\ref{prop:interior} come from upper estimates of the expected loss at a uniform threshold. The rounding rule returns the least loss over all its candidate subspaces. For $v_*>0$, how close $V(T)/v_*$ can come to these bounds at an optimal relaxed solution remains open here.
	\end{remark}
	
	\section{The bound for arbitrary symmetric inputs}\label{sec:arbitrary}
	Fix a finite feasible triple $(T_0,U,W)$. We first vary the eigenvalues of $T_0$ in a fixed basis, keeping $U,W$ fixed. If the eigenvalues stay in order and each original eigenspace stays intact, every threshold projector of the new matrix already belongs to $\mathcal P(T_0)$. We use this freedom to minimize $L$ over the admissible eigenvalues. The derivative at the minimizer gives an identity that cancels two terms in the expected loss estimate.
	
	\subsection{Stationarity with respect to the eigenvalues}
	Choose an eigenbasis of $T_0$ in decreasing order of eigenvalue. Minimize $L(T,U,W)$ over matrices $T$ diagonal in this basis, keeping $U,W$ fixed. Require the diagonal entries of $T$ to lie in $[0,1]$ in nonincreasing order, with equal entries wherever the corresponding eigenvalues of $T_0$ are equal. This compact set contains $T_0$, so lower semicontinuity gives a minimizer $\bar T$ with finite loss. Adjacent distinct eigenvalues may become equal at $\bar T$. A threshold then includes both original eigenspaces or neither of them. Hence $\mathcal P(\bar T)\subseteq\mathcal P(T_0)$ and
	\begin{equation}\label{eq:order}
		V(T_0)\le V(\bar T),\qquad L(\bar T,U,W)\le L(T_0,U,W).
	\end{equation}
	Use the chosen eigenbasis throughout this section, and write $x_i$ for the eigenvalues of $\bar T$. To differentiate along an admissible perturbation, we use the polynomial
	\begin{align*}
		\lambda(x)&=x(1-x)(2x-1)
		\big[-0.001447+0.213196(2x-1)^2-0.166343(2x-1)^4\big].
	\end{align*}
	Throughout, terminating decimals used as coefficients denote exact rational numbers. The perturbation $x\mapsto x+\epsilon\lambda(x)$ fixes both endpoints. Since $\lambda'$ is bounded, the map is strictly increasing for all sufficiently small $|\epsilon|$. Thus an interior eigenvalue stays in $(0,1)$ and the eigenvalues retain their order and ties. By~\eqref{eq:k}, zero multipliers occur only at pairs of endpoints, so the positions of all zero multipliers $k(x_i,x_j)$ are preserved. Write $\partial_1 k,\partial_2 k$ for the partial derivatives with respect to the first and second arguments. For $k(x,y)>0$ and $(x,y)\ne(1,1)$, define
	\[
	J_\lambda(x,y)=
	\frac{\lambda(x)\partial_1 k(x,y)+\lambda(y)\partial_2 k(x,y)}{k(x,y)^2}.
	\]
	Set $J_\lambda(1,1)=0$ because the perturbation fixes $1$. Set $J_\lambda(x,y)=0$ also when $k(x,y)=0$, solely to include those indices in the sums below. The corresponding matrix entries vanish by finiteness of the quadratic form, and the zero multipliers stay zero along the perturbation. For $k>0$, the derivative of $1/k$ along the perturbation at $\epsilon=0$ is $-J_\lambda(x,y)$. Also, $1-(x+\epsilon\lambda(x))=(1-x)+\epsilon\lambda(1-x)$ because $\lambda(1-x)=-\lambda(x)$. Both signs of $\epsilon$ are feasible at the minimizer $\bar T$. The derivative of $L$ along this curve at $\epsilon=0$ therefore vanishes. Only the two quadratic terms vary, giving
	\begin{equation}\label{eq:stationary}
		\sum_{ij}J_\lambda(x_i,x_j)U_{ij}^2+
		\sum_{ij}J_\lambda(1-x_i,1-x_j)W_{ij}^2=0 .
	\end{equation}
	Set $c=1.189$ for the remainder of this section. By~\eqref{eq:order}, it suffices to prove $V(\bar T)\le cL(\bar T,U,W)$.
	
	\subsection{Expected loss after thresholding and diagonal scaling}
	Choose $s$ uniformly on $[0,1]$ and put $P_s=P_s(\bar T)$ in this section. We construct a fit for $C$ by scaling $U$ on $\ran P_s$. For $x,s\in[0,1]$, define
	\begin{align*}
		p(x)&=0.743787-0.486537x+0.627925x^2-0.392282x^3,\\
		g(x)&=(0.426605-2.315288x+2.696582x^2)/(0.12+x),\\
		\rho_s(x)&=\frac{1+(1-x)g(x)(x-2s)}{1-(1-x)p(x)} .
	\end{align*}
	Let $R_s$ be diagonal in the chosen eigenbasis. Its entry at eigenvalue $x$ is $\rho_s(x)$ if $x\ge s$ and zero otherwise. The denominator of $\rho_s(x)$ is positive because $p$ decreases on $[0,1]$ and $p(0)<1$. Then $R_sUR_s$ is positive semidefinite and has range in $\ran P_s$. For $D$, apply the same rule to $W$ using $I-\bar T$ and threshold $1-s$. Unless $s$ is an eigenvalue of $\bar T$, the condition $1-x\ge1-s$ is equivalent to $x<s$. The subspace selected for $W$ is therefore $\ran(I-P_s)$ almost surely.
	
	For the entry $(i,j)$ with $x_i,x_j>0$, put $x=\max(x_i,x_j)$ and $y=\min(x_i,x_j)$. For this pair let $h_s=\rho_s(x)\rho_s(y)$. The $(i,j)$ entry of $R_sUR_s$ is $h_s U_{ij}$ for $s\le y$ and zero otherwise. Its contribution to $\E\,\fn{C-R_sUR_s}^2$ is
	\[
	C_{ij}^2-2C_{ij}U_{ij}\int_0^y h_s\,ds+U_{ij}^2\int_0^y h_s^2\,ds.
	\]
	Subtract this expected loss from the scaled relaxed loss $c(C_{ij}^2-2C_{ij}U_{ij}+U_{ij}^2/k(x,y))$. For fixed $U_{ij}$, the minimum of the difference over $C_{ij}\in\mathbb R$ is $M(x,y)U_{ij}^2$, where
	\begin{equation}\label{eq:M}
		M(x,y)=\frac{c}{k(x,y)}-\int_0^y h_s^2\,ds
		-\frac{(c-\int_0^y h_s\,ds)^2}{c-1}.
	\end{equation}
	For almost every $s$, the two scaled matrices are feasible on $\ran P_s$ and $\ran(I-P_s)$. Their loss therefore bounds $v(P_s)$ from above by \eqref{eq:hard}.
	
	\begin{lemma}[A pointwise lower bound for $M+J_\lambda$]\label{lem:gen}
		Define $f$ on $[0,1]$ and $\theta$ on $[0,1)$ by
		\begin{align*}
			f(x)&=(-0.398243+0.543957x-0.533234x^2)/(0.057958+x),\\
			\theta(x)&=2.191072\log(1-x)-1.031425x-0.481436x^2 .
		\end{align*}
		For $x,y\in[0,1)$, put
		\[
		 G(x,y)=\sqrt{(1-x)(1-y)}\,f(x)f(y)\cos\bigl(\theta(x)-\theta(y)\bigr).
		\]
		Set $G(x,y)=0$ if $x=1$ or $y=1$. Then for $0<y\le x\le1$,
		\begin{equation}\label{eq:minorant}
			M(x,y)+J_\lambda(x,y)\ge G(x,y).
		\end{equation}
	\end{lemma}
	\begin{proof}
		For $y<1$, set $b=1-y$ and $t=\sqrt{(1-x)/(1-y)}$. Since $b,k(x,y)>0$, the desired inequality is equivalent to $(k(x,y)/b)(M+J_\lambda-G)\ge0$. Appendix~\ref{app:general} gives a positive lower bound for this scaled difference on each rectangle in a cover of $[0,1]^2$. At $x=y=1$ we have $k=h_s=1$, so \eqref{eq:M} gives $M=0$. By definition $J_\lambda=G=0$ there as well.
	\end{proof}
	
	\begin{proof}[Proof of Theorem~\ref{thm:main} for arbitrary inputs]
		The matrix $[G(x_i,x_j)]$ is a Gram matrix. Take the zero vector at $x=1$ and the vector $\sqrt{1-x}\,f(x)(\cos\theta(x),\sin\theta(x))$ at $x<1$. The inner product of the vectors at $x_i,x_j$ equals $G(x_i,x_j)$. For matrices of equal size let $\circ$ denote the entrywise product. By the Schur product theorem~\cite[Sec.~7.5]{HornJohnson2013}, $U\circ U=[U_{ij}^2]$ is also positive semidefinite. To see this directly, let $u_i^\top$ be row $i$ of $U^{1/2}$. Then $U_{ij}=u_i^\top u_j$ and $U_{ij}^2=\ip{u_i u_i^\top}{u_j u_j^\top}$. The trace inner product of two positive semidefinite matrices $A,B$ is nonnegative, since $\tr(AB)=\tr(A^{1/2}BA^{1/2})\ge0$. Hence $\sum_{ij}G(x_i,x_j)U_{ij}^2\ge0$.

		If $x_i=0$, finiteness of $\qf_{\bar T}(U)$ forces $U_{ii}=0$. The inequality $U_{ij}^2\le U_{ii}U_{jj}$ then gives $U_{ij}=0$ for every $j$. Each entry in that row or column contributes $(c-1)C_{ij}^2\ge0$ to $c\Phi_C(\bar T,U)-\E\,\fn{C-R_sUR_s}^2$. For entries with both eigenvalues positive, multiply \eqref{eq:minorant} by $U_{ij}^2$ and sum. Since the sum involving $G$ is nonnegative, the loss difference for $C$ is at least $-\sum_{ij}J_\lambda(x_i,x_j)U_{ij}^2$. The corresponding lower bound for $D$ uses $W$ and the eigenvalues $1-x_i$. The two sums involving $J_\lambda$ cancel by \eqref{eq:stationary}. Thus the expected loss of the constructed pair is at most $cL(\bar T,U,W)$. Since $V(\bar T)\le\E\,v(P_s)$, \eqref{eq:order} gives $V(T_0)\le V(\bar T)\le cL(\bar T,U,W)\le cL(T_0,U,W)$.
	\end{proof}
	
	\section{Positive semidefinite inputs after reflection}\label{sec:positive}
	Fix a finite feasible triple $(T,U,W)$ and suppose $C,D\succeq0$. Choose $s$ uniformly on $[0,1]$ and write $P_s=P_s(T)$. The compressions $P_s C P_s$ and $(I-P_s)D(I-P_s)$ are positive semidefinite and are therefore the minimizing pair in \eqref{eq:hard}.
	
	Fix an eigenbasis of $T$, write $x_i$ for its eigenvalues, and put $a_i=1-x_i$. The compression $P_sCP_s$ retains entry $C_{ij}$ precisely when $s\le\min(x_i,x_j)$. The expected squared loss from this entry is therefore $\max(a_i,a_j)C_{ij}^2$. For $k=k(x_i,x_j)>0$, completing the square gives the corresponding term in $\Phi_C(T,U)$ as
	\[
	C_{ij}^2-2C_{ij}U_{ij}+U_{ij}^2/k
	=(1-k)C_{ij}^2+(U_{ij}-kC_{ij})^2/k .
	\]
	If $k=0$, finiteness of $\qf_T(U)$ requires $U_{ij}=0$. In both cases the contribution to the relaxed loss is at least $(1-k)C_{ij}^2$. Equation~\eqref{eq:k} gives $1-k=\beta(a_i,a_j)$, where
	\[
	\beta(a,b)=\frac{a^2+b^2}{a+b}\quad(a+b>0),\qquad \beta(0,0)=0.
	\]
	For $a,b>0$, put $t=|\log(a/b)|$ and $c_+=1.15487$. Dividing $c_+\beta(a,b)-\max(a,b)$ by $\sqrt{ab}$ yields a function of $t$ alone,
	\begin{equation}\label{eq:F}
		\frac{c_+\beta(a,b)-\max(a,b)}{\sqrt{ab}}
		=F(t):=c_+\frac{\cosh t}{\cosh(t/2)}-e^{t/2}.
	\end{equation}
	
	\begin{lemma}[A lower bound by a cosine sum]\label{lem:cos}
		Let $(\omega_m,d_m)$, $m=1,\ldots,98$, be the nonnegative rational frequency and weight pairs in Table~\ref{tab:cosine-data}. Then
		\begin{equation}\label{eq:cos}
			F(t)\ge \sum_{m=1}^{98} d_m\cos(\omega_m t)\qquad (t\ge0).
		\end{equation}
	\end{lemma}
	\begin{proof}
		Let $\delta(t)=F(t)-\sum_m d_m\cos(\omega_m t)$. Appendix~\ref{app:cosine} gives $\delta(j/4000)\ge2535846/10^{12}$ for $j=0,\ldots,12000$ and $|\delta''(t)|\le54$ on $[0,3]$. The linear interpolant of $\delta$ between adjacent nodes is at least the smaller endpoint value. On a node interval $[a,b]$, the linear interpolation remainder is $\delta''(\xi)(t-a)(t-b)/2$ for some $\xi\in(a,b)$~\cite[eq.~(3.3.5)]{TemmeDLMF}. Its magnitude is at most $54/(8\cdot4000^2)$, less than the lower bound for $\delta$ at the nodes. Thus $\delta\ge0$ on $[0,3]$. For $t\ge3$, Appendix~\ref{app:cosine} gives $F(t)>\sum_m d_m\ge\sum_m d_m\cos(\omega_m t)$.
	\end{proof}
	\begin{proof}[Proof of Theorem~\ref{thm:main} when $C,D\succeq0$]
		For each $a_i>0$, form a vector by concatenating the coordinate pairs
		\[
		\sqrt{a_i d_m}\big(\cos(\omega_m\log a_i),\sin(\omega_m\log a_i)\big)
		\]
		over the frequencies $m$. For $a_i=0$ take the zero vector. Let $K$ be the Gram matrix of these vectors. If $a_i,a_j>0$, the cosine difference identity and \eqref{eq:F}--\eqref{eq:cos} give $K_{ij}\le c_+\beta(a_i,a_j)-\max(a_i,a_j)$. If $a_i$ or $a_j$ is zero, then $K_{ij}=0$ and the right side equals $(c_+-1)\max(a_i,a_j)\ge0$.

		Multiply the bound on $K_{ij}$ by $C_{ij}^2$ and sum. As $C\succeq0$, the argument for $U\circ U$ in Section~\ref{sec:arbitrary} gives $C\circ C\succeq0$. Its trace inner product with $K$ is nonnegative. It follows that $\sum_{ij}\max(a_i,a_j)C_{ij}^2\le c_+\sum_{ij}\beta(a_i,a_j)C_{ij}^2$. The lower bound for $\Phi_C(T,U)$ now gives
		\[
		\E\,\fn{C-P_s C P_s}^2\le c_+\Phi_C(T,U).
		\]
		Apply the same estimate to $D,W$ using $I-T$ and threshold $1-s$. The corresponding projector is $I-P_s$ almost surely. Summing the two estimates gives $\E\,v(P_s)\le c_+ L(T,U,W)$, and $V(T)\le\E\,v(P_s)$ by \eqref{eq:V}.
	\end{proof}
	
	\section{The bound for at most two distinct interior eigenvalues}\label{sec:twolevel}
	For two distinct eigenvalues in $(0,1)$, we use a separate scaling on each interior eigenspace. Positive semidefiniteness then bounds the Frobenius norm of the block between the two eigenspaces in terms of the norms of the diagonal blocks.
	
	\begin{proposition}\label{prop:interior}
		Suppose that $T$ has at most two distinct eigenvalues in $(0,1)$. For every finite feasible triple $(T,U,W)$ one has $V(T)\le(59/50)L(T,U,W)$. Eigenvalues $0$ and $1$ may both be present, and their multiplicities and those of the interior eigenvalues are unrestricted.
	\end{proposition}
	\begin{proof}
		Set $c=59/50$ and $\alpha=c-1$. First suppose there are two distinct interior eigenvalues $0<y<x<1$, and let $r=y/x$. Choose $\sigma$ uniformly on $[0,1]$ and write $P_\sigma=P_\sigma(T)$. For a parameter $\kappa\in\{0,\tfrac12,1\}$, put
		\[
		u=1+\kappa(1-r)(1-x),\qquad
		b=1-r\kappa(1-x)[2+\kappa(1-r)(1-x)].
		\]
		For $b\ge0$, define a diagonal matrix $R$ in an eigenbasis of $T$. On the $y$-eigenspace its diagonal entry is $1/\sqrt y$ if $\sigma\le y$ and zero otherwise. On the $x$-eigenspace the entry is $u/\sqrt x$ for $\sigma\le y$, $\sqrt{b/x}$ for $y<\sigma\le x$, and zero for $\sigma>x$. On the eigenspaces at $0$ and $1$ use zero and one, respectively, for every $\sigma$. Then $RUR$ is positive semidefinite with range in $\ran P_\sigma$. Finiteness of $\qf_T(U)$ and positive semidefiniteness of $U$ force its rows and columns at eigenvalue $0$ to vanish.
		
		Fix $\xi,\eta\in\{x,y,1\}$ with $k(\xi,\eta)>0$. Let $h_\sigma$ be the product of the diagonal entries of $R$ assigned to $\xi$ and $\eta$ at threshold $\sigma$. Define
		\[
		m_{\xi\eta}=\frac{c}{k(\xi,\eta)}-\E h_\sigma^2-\frac{(c-\E h_\sigma)^2}{\alpha}.
		\]
		For an entry $U_{ij}$ whose row and column correspond to eigenvalues $\xi$ and $\eta$, minimization over $C_{ij}$ as in \eqref{eq:M} gives $m_{\xi\eta}U_{ij}^2$ as a lower bound on the entry's contribution to $c\Phi_C(T,U)-\E\,\fn{C-RUR}^2$. Appendix~\ref{app:twolevel} uses Bernstein coefficients to show that for every pair $(x,y)$ one can choose $\kappa$ so that $b\ge0$ and
		\begin{equation}\label{eq:two_scalar}
			m_{xx},m_{yy},m_{1x},m_{1y}\ge0,\qquad
			m_{xy}+\sqrt{m_{xx}m_{yy}}\ge0 .
		\end{equation}
		Entries involving a zero eigenvalue contribute $(c-1)C_{ij}^2\ge0$ to the loss difference. For the eigenspace at $1$, both multipliers are one, so $m_{11}=0$. The bounds on $m_{1x}$ and $m_{1y}$ in \eqref{eq:two_scalar} cover the other entries involving $1$. For $\xi,\eta\in\{x,y\}$ let $U_{\xi\eta}$ denote the block between those eigenspaces. The two off-diagonal blocks are transposes, so their combined contribution is $2m_{xy}\fn{U_{xy}}^2$. When $m_{xy}<0$, factor $U=\Gamma\Gamma^\top$ and let $\Gamma_x,\Gamma_y$ be the corresponding row blocks. Then $\fn{U_{xy}}^2=\ip{\Gamma_x^\top\Gamma_x}{\Gamma_y^\top\Gamma_y}$, so Cauchy--Schwarz gives $\fn{U_{xy}}^2\le\fn{U_{xx}}\fn{U_{yy}}$. By \eqref{eq:two_scalar}, the total contribution from the four blocks on the $x$- and $y$-eigenspaces is at least
		\[
		\left(\sqrt{m_{xx}}\,\fn{U_{xx}}-\sqrt{m_{yy}}\,\fn{U_{yy}}\right)^2 .
		\]
		If $m_{xy}\ge0$, all four block contributions are nonnegative. If $x$ is the only interior eigenvalue, choose any $y\in(0,x)$ and take the auxiliary $y$-subspace to be zero-dimensional. Every term involving that subspace vanishes. If $T=P$ is a projector, then $P\in\mathcal P(P)$, and finiteness of $\qf_P(U)$ and $\qf_{I-P}(W)$ gives $U=PUP$ and $W=(I-P)W(I-P)$. Hence $V(P)\le v(P)\le L(P,U,W)$.

		For $D,W$, apply the construction to $I-T$ with threshold $1-\sigma$. The two scalings may use different values of $\kappa$ because $\kappa$ affects only the multipliers within the selected subspaces. These subspaces are complementary almost surely. Adding the two loss bounds gives $V(T)\le\E\,v(P_\sigma)\le(59/50)L(T,U,W)$.
	\end{proof}
	For $n\le2$, $T$ automatically has at most two distinct eigenvalues in $(0,1)$. In higher dimensions the spectral constraint $0\preceq T\preceq I$ allows more. For arbitrary symmetric inputs it remains open whether $V(T_*)\le1.18v_*$ for every $T_*$ occurring in an optimal solution of \eqref{eq:relax}.
	
	\section*{Acknowledgments}
	At the author's research talk on 21 November 2026, Houduo Qi (Department of Applied Mathematics, The Hong Kong Polytechnic University) asked whether an approximation algorithm was available for this projection problem. The author thanks him for the question, which led him to write this note later that evening. The author used generative AI in preparing the note and takes full responsibility for its content.
	
	\appendix
	\section{Verification of the scalar inequalities}\label{app:verification}
	
	\subsection{Interval bounds for the pointwise inequality}\label{app:general}
	Fix $c=1.189$, with $p,g,\lambda$ as in Section~\ref{sec:arbitrary} and $f,\theta$ as in Lemma~\ref{lem:gen}. For $0<y<1$ and $y\le x\le1$, set $b=1-y$, $t=\sqrt{(1-x)/(1-y)}$, and $u=t^2$. Then $y=1-b$, $x=1-bu$, $0<b<1$, and $0\le t\le1$. We evaluate interval bounds on the closed square $[0,1]^2$ in $(b,t)$. In this calculation $k,M,J_\lambda,G$ without arguments denote their values at $(x,y)$. With $v=y+ux$, we have
	\[
		k=\frac{v}{1+u}=1-b\,\frac{1+u^2}{1+u}.
	\]
	For the product $h_s$ in \eqref{eq:M}, write $h_s=1+bH_s$. With $a_s=u g(x)(x-2s)$, $b_s=g(y)(y-2s)$, and $D_0=(1-bu p(x))(1-bp(y))$, we have
	\[
		H_s=\frac{u p(x)+p(y)(1-bu p(x))}{D_0}
		   +\frac{a_s+b_s(1+b a_s)}{D_0}.
	\]
	The polynomial $p$ decreases on $[0,1]$ because its derivative is a quadratic with negative leading coefficient and negative discriminant. Since $p(0)<1$, both factors of $D_0$ are positive. Use the three-point Gauss--Legendre rule on $[0,1]$~\cite[Sec.~3.5(v)]{TemmeDLMF}, with nodes $(r_1,r_2,r_3)=(\tfrac12-\sqrt{3/20},\,\tfrac12,\,\tfrac12+\sqrt{3/20})$ and weights $(w_1,w_2,w_3)=(5/18,\,4/9,\,5/18)$. For $i=1,2$ put $m_i=\sum_{j=1}^3 w_j(H_{y r_j})^i$. The moment identities $\sum_j w_j r_j^\ell=1/(\ell+1)$ for $\ell=0,\ldots,4$ give $\int_0^y H_s^i\,ds=y m_i$, since $H_s$ is quadratic in $s$.

	For $t>0$ take $\omega=2.191072$ and use the phase
	\[
		-2\omega\log t+(y-x)\big[-1.031425-0.481436(x+y)\big].
	\]
	Define $G_0$ to be $t f(x)f(y)$ times the cosine of this phase for $t>0$, and set $G_0=0$ at $t=0$. When $b,t>0$, the phase equals $\theta(y)-\theta(x)$, because $(1-x)/(1-y)=t^2$. Hence $G=bG_0$. The explicit phase also defines $G_0$ on $b=0$ for $t>0$ without evaluating $\theta(1)$.

	Put $C_0=y m_2+(y m_1-1)^2/(c-1)$ and, for $b,k>0$, put $D_\lambda=(k/b)J_\lambda$. Substituting $h_s=1+bH_s$ into \eqref{eq:M} and collecting powers of $b$ gives
	\begin{equation}\label{eq:interval-target}
		\frac{k}{b}(M+J_\lambda-G)=c\,\frac{1+u^2}{1+u}-k-k(b C_0+G_0)+D_\lambda.
	\end{equation}
	Let $\ell$ be the polynomial obtained by removing the factor $x(1-x)$ from $\lambda(x)$, so that $\ell$ is defined also at $0$ and $1$. For the derivative term put
	\[
		d_x=\ell(x)\frac{u^2+2u-1}{1+u},\qquad
		d_y=\ell(y)\frac{1+2u-u^2}{1+u}.
	\]
	Differentiating $k$ gives $D_\lambda=(ux\,d_x+y\,d_y)/v$ wherever $b,k>0$. Use this expression wherever $v>0$. Its value is then a convex combination of $d_x$ and $d_y$. On each rectangle the smaller of the lower bounds for $d_x,d_y$ therefore bounds $D_\lambda$ from below at every point with $v>0$. This remains valid when the interval for $v$ contains zero and division by that interval is unavailable. Equation~\eqref{eq:interval-target} is needed only for $b,k>0$. The case $b=0$ gives $x=y=1$, treated directly in Lemma~\ref{lem:gen}.

	Let $[t_0,t_1]$ be the $t$-interval of a rectangle. If $t_0<0.28$, bound the cosine in $G_0$ by $[-1,1]$. Otherwise $t_0>0$ and we bound the logarithm in the explicit phase as follows.
	At the midpoint $m$ of $[t_0,t_1]$, put $z=(1-m)/(1+m)\le 9/16$. The series for $\log((1+z)/(1-z))=-\log m$ gives the approximation $4\omega\sum_{j=0}^{55}z^{2j+1}/(2j+1)$ to $-2\omega\log m$. The omitted tail is bounded by $4\omega z^{113}/[113(1-z^2)]<10^{-16}$. For $t\in[t_0,t_1]$, the mean value theorem gives $|2\omega(\log t-\log m)|\le\omega(t_1-t_0)/t_0$. Adding this bound to the truncation error encloses $-2\omega\log t$ throughout the interval.

	For the cosine, let $a$ be the nearest multiple of $2^{-15}$ to the midpoint of the interval for its argument. These values satisfy $|a|<10$. Starting with the constant term of the Maclaurin series, multiply the preceding term by $-a^2/[(2j)(2j-1)]$ at step $j$, beginning with $j=1$. Sum through degree $46$ and enclose the remainder by $10^{48}/48!<10^{-12}$ in absolute value. To bound the cosine on the whole argument interval, enlarge the enclosure for $\cos a$ by the distance from $a$ to the farther endpoint, using the $1$-Lipschitz property. Intersect the result with $[-1,1]$.

	Use interval arithmetic for sums, differences and products so that each computed interval contains all attainable values~\cite[Chs.~2--4]{MooreKearfottCloud2009}. When squaring an interval containing zero, take zero as the lower endpoint. Divide by a positive interval by multiplying by its reciprocal. Evaluate the coefficient polynomials and truncated logarithm series by Horner's rule. Adjust each computed binary64 endpoint $a$ outward by $(1+|a|)2^{-49}$ to allow for rounding of the endpoint and the adjustment itself. Use the same allowance for the input constants and Gauss nodes. Squaring shows that the error in approximating $\sqrt{3/20}$ by $0.3872983346207417$ is below $10^{-16}$.

	Start with the square $[0,1]^2$. On each rectangle compute $k=v/(1+u)$. If the interval for $D_0$ has a nonpositive lower endpoint, bisect the rectangle. Otherwise form $1/D_0$ once for the two terms of $H_s$. For $D_\lambda$, first use the smaller of the lower endpoints for $d_x,d_y$. If the lower endpoint for $v$ is positive, evaluate $(ux\,d_x+y\,d_y)/v$ as well and keep the larger of the two lower bounds. Accept the rectangle if the lower bound on the right of \eqref{eq:interval-target} is strictly positive. Otherwise bisect it. Each bisection is at the midpoint, alternating $b$ and $t$ and starting with $b$. The $11\,839\,610$ accepted rectangles cover the square. Their maximum depth is $31$, and the smallest lower bound exceeds $1.7\times10^{-10}$.

	\subsection{Verification of the cosine lower bound}\label{app:cosine}
	Let $\delta(t)=F(t)-\sum_m d_m\cos(\omega_m t)$, with $F$ given by \eqref{eq:F} and $c_+=1.15487$. The frequencies and weights in Table~\ref{tab:cosine-data} give $\mu:=\sum_m d_m<4/25$ and $\sum_m d_m\omega_m^2<50$. The total weight $\mu$ bounds the cosine sum from above for $t\ge3$. The second sum bounds the absolute second derivative of the cosine sum on $[0,3]$. We use it to extend a lower bound on $\delta$ at the nodes to the intervals between them.
	\begin{table}[!htbp]
		\centering
		\caption{Data for the cosine bound \eqref{eq:cos}. Each adjacent pair of columns gives the scaled frequency and weight shown in its heading. Read the pairs from left to right in each row.}
		\label{tab:cosine-data}
		\small\setlength{\tabcolsep}{3.8pt}\renewcommand{\arraystretch}{0.94}
		\begin{tabular}{@{}rr@{\quad}rr@{\quad}rr@{\quad}rr@{\quad}rr@{}}
			\toprule
			$100\omega_m$ & $10^{10}d_m$ & $100\omega_m$ & $10^{10}d_m$ & $100\omega_m$ & $10^{10}d_m$ & $100\omega_m$ & $10^{10}d_m$ & $100\omega_m$ & $10^{10}d_m$\\
			\midrule
		297 & 132530224 & 312 & 727610379 & 578 & 53987891 & 593 & 204854308 & 866 & 32656035 \\
		881 & 84909950 & 896 & 659902 & 1156 & 45136331 & 1186 & 22232023 & 1447 & 29884142 \\
		1477 & 13113541 & 1737 & 20568717 & 1767 & 9213989 & 2028 & 15637737 & 2058 & 6156926 \\
		2320 & 12938571 & 2350 & 3655984 & 2611 & 10571697 & 2641 & 2470659 & 2902 & 8812243 \\
		2932 & 1688686 & 3193 & 7439563 & 3223 & 1184365 & 3485 & 6592322 & 3515 & 596429 \\
		3776 & 5682429 & 3806 & 399829 & 4067 & 4904453 & 4097 & 301871 & 4359 & 4391955 \\
		4389 & 100622 & 4650 & 3810937 & 4680 & 104577 & 4942 & 3422173 & 4972 & 10216 \\
		5233 & 2983613 & 5263 & 50448 & 5525 & 2690875 & 5816 & 2349846 & 5831 & 51934 \\
		6107 & 2051889 & 6122 & 102072 & 6399 & 1928833 & 6429 & 8925 & 6690 & 1719182 \\
		6720 & 32956 & 6982 & 1585608 & 7012 & 1771 & 7273 & 1418523 & 7303 & 26625 \\
		7565 & 1317605 & 7856 & 1183139 & 7886 & 21876 & 8148 & 1101216 & 8439 & 995000 \\
		8469 & 15422 & 8731 & 925759 & 9022 & 837907 & 9052 & 11290 & 9314 & 778576 \\
		9605 & 700447 & 9635 & 15101 & 9897 & 653642 & 9927 & 3180 & 10188 & 588757 \\
		10218 & 15274 & 10480 & 553920 & 10771 & 505787 & 10801 & 3121 & 11063 & 466742 \\
		11354 & 411623 & 11369 & 14789 & 11646 & 389134 & 11937 & 343807 & 11952 & 11385 \\
		12228 & 314534 & 12258 & 8956 & 12520 & 291648 & 12811 & 258416 & 12841 & 4092 \\
		13103 & 235801 & 13394 & 200668 & 13409 & 9619 & 13686 & 185786 & 13977 & 157653 \\
		14007 & 3685 & 14269 & 138679 & 14560 & 109888 & 14575 & 6376 & 14852 & 95887 \\
		15143 & 74052 & 15173 & 2461 & 15435 & 54110 & 15465 & 3257 & 15726 & 27727 \\
		15756 & 10203 & 16018 & 11849 & 16048 & 7260 &  &  &  &  \\
			\bottomrule
		\end{tabular}
	\end{table}

\begin{samepage}
	On $[0,3]$, set $N=4000$ and use the nodes $t_j=j/N$ for $j=0,\ldots,3N$. Write $\mathrm i$ for the imaginary unit and $\operatorname{Re}$ for the real part. Put $\zeta_m=\omega_m/N$. At successive nodes $e^{\mathrm i\omega_m t_j}$ is multiplied by $e^{\mathrm i\zeta_m}$. To approximate that multiplier, fix $S=10^{24}$ and let $R_S$ round each of the real and imaginary parts down to a multiple of $1/S$. Define
	\[
		q_m=R_S\!\left(\sum_{\ell=0}^{17}\frac{(\mathrm i\zeta_m)^\ell}{\ell!}\right),
		\qquad z_{m,0}=1,\qquad z_{m,j+1}=R_S(q_m z_{m,j}) .
	\]
\end{samepage}
	The exponential remainder is at most $\zeta_m^{18}/18!<1/S$, and the error from rounding both coordinates is below $2/S$. Hence $|q_m-e^{\mathrm i\zeta_m}|<3/S$. Suppose $|z_{m,j}-e^{\mathrm i j\zeta_m}|\le8j/S$. Since $j\le3N$ and $8(3N)/S<1$, we have $|z_{m,j}|\le2$. At the next step the multiplier error contributes at most $2(3/S)$ and rounding contributes at most $2/S$. Starting from $z_{m,0}=1$, this proves $|z_{m,j}-e^{\mathrm i j\zeta_m}|\le8j/S$ for all $j\le3N$. Taking real parts and summing with weights $d_m$ gives the upper bound $\sum_m d_m\operatorname{Re}(z_{m,j})+8j\mu/S$ for the cosine sum at node $j$.
	
\begin{samepage}
	Let $L_j$ be $\sum_{\ell=0}^{32}(t_j/2)^\ell/\ell!$ rounded down to a multiple of $1/S$, and put $U_j=L_j+2/S$. On this interval the ratio of successive omitted terms is less than one half, so the positive remainder is bounded by $2(3/2)^{33}/33!<1/S$. Thus $L_j\le e^{t_j/2}\le U_j$. Substitution of $z=e^{t/2}$ into \eqref{eq:F} gives
	\begin{equation}\label{eq:F-exp}
		F(t)=(c_+-1)z+\frac{c_+}{z}-\frac{2c_+z}{z^2+1}
	\end{equation}
\end{samepage}
	For $z\ge1$ the fraction $z/(z^2+1)$ decreases. Replacing $z$ by $L_j$ in the last term of \eqref{eq:F-exp} therefore bounds that term from below. Subtracting the upper bound on the cosine sum gives a rational lower bound for $\delta(t_j)$. For $j=0,\ldots,3N$, this lower bound satisfies
	\[
		(c_+-1)L_j+\frac{c_+}{U_j}-\frac{2c_+L_j}{L_j^2+1}
		-\sum_m d_m\operatorname{Re}z_{m,j}-\frac{8j\mu}{S}
		\ \ge\ \frac{2535846}{10^{12}} .
	\]
	Write $\operatorname{sech}u=1/\cosh u$. On $[0,3]$ we have $\cosh(t/2)<3$, $e^{t/2}<5$, and $|\operatorname{sech}''(t/2)|\le1$, where the primes denote derivatives with respect to the argument of $\operatorname{sech}$. Differentiating $\cosh t/\cosh(t/2)=2\cosh(t/2)-\operatorname{sech}(t/2)$ twice gives an absolute bound of $7/4$ on the second derivative of the ratio. The second derivative of $e^{t/2}$ is at most $5/4$. Together with $\sum_m d_m\omega_m^2<50$, these bounds give $|\delta''|<54$, as needed for the interpolation estimate in Lemma~\ref{lem:cos}.

	For $t\ge3$, $z=e^{t/2}>4$. Equation~\eqref{eq:F-exp} gives $F(t)\ge(c_+-1)z-c_+/z\ge4(c_+-1)-c_+/4>4/25>\mu$. Since $\sum_m d_m\cos(\omega_m t)\le\mu$, this proves \eqref{eq:cos} on $[3,\infty)$.
	
	\subsection{Polynomial bounds for two interior eigenvalues}\label{app:twolevel}
	Fix $0<y<x<1$ and $\kappa\in\{0,\tfrac12,1\}$. Set $c=59/50$, $\alpha=c-1$, $r=y/x$, and $s=\sqrt r$, and take $u,b$ as in the proof of Proposition~\ref{prop:interior}. The coefficients $m_{\xi\eta}$ involve $\sqrt b$ and require $b\ge0$. The polynomial tests below are defined for every candidate $\kappa$ and include $b\ge0$ among the conditions to be checked.

	For $b\ge0$, multiplying the three interior coefficients by $x$ gives
	\begin{align*}
		d=xm_{xx}&=\alpha(1-x)-r(1-r)\kappa^2(1-x)^2[2+\kappa(1-r)(1-x)]^2,\\
		xm_{yy}&=\alpha(1-y)/r,\\
		xm_{xy}&=cx/k(x,y)-u^2-x(c-us)^2/\alpha .
	\end{align*}
	Put $n_0=1-x+r(1-y)>0$. Then $x/k(x,y)=(2-x-y)/n_0$. To check the mixed coefficient, form the polynomials
	\[
	 \chi=\alpha c(2-x-y)-n_0[\alpha u^2+x(c-us)^2],
	 \qquad \Delta=\alpha(1-y)d(\alpha n_0)^2-r\chi^2 .
	\]
	Since $\chi=\alpha n_0 x m_{xy}$, it has the sign of $m_{xy}$. Once $d\ge0$, both $m_{xx}$ and $m_{yy}$ are nonnegative. If $\chi\ge0$, $m_{xy}+\sqrt{m_{xx}m_{yy}}\ge0$ follows immediately. If $\chi<0$, the inequality $\Delta\ge0$ suffices because $\Delta=r(\alpha n_0 x)^2(m_{xx}m_{yy}-m_{xy}^2)$. For the coefficient $m_{1x}$ put
	\[
	m_0=1-\tfrac12\{r(u-1)^2+(1-r)(b-1)^2\}.
	\]
	Conditional on $\sigma\le x$, the event $\sigma\le y$ has probability $r=y/x$. Thus the diagonal entry of $\sqrt{x}\,R$ on the $x$-eigenspace equals $u$ or $\sqrt b$ with probabilities $r$ and $1-r$, respectively. Direct expansion gives $ru^2+(1-r)b=1$. Expanding the squared deviations from one and using $(\sqrt b-1)^2\le(b-1)^2$ gives the lower bound in $m_0\le ru+(1-r)\sqrt b\le1$. The upper bound follows from Cauchy--Schwarz. In the definition of $m_{1x}$ we have $k(1,x)=x$, $\E h_\sigma^2=1$, and $\E h_\sigma=\sqrt{x}\,[ru+(1-r)\sqrt b]$. Define $E=\alpha(c-x)-x(c-\sqrt{x}\,m_0)^2$. Then $E\ge0$ suffices for $m_{1x}\ge0$. For $m_{1y}$ evaluate $E$ at $r=1$ with $x$ replaced by $y$. Then $m_0=1$ and this bound is independent of $\kappa$.

	Substituting $x=t^2$ with $t=\sqrt{x}$ and $r=s^2$ makes $b,d,\chi,\Delta,E$ polynomials in $(t,s)$. Parametrize the triangle $t\ge s$ by $t=1-(1-s)z$ and the triangle $s\ge t$ by $s=1-(1-t)z$, with $z\in[0,1]$ in both cases. Divide out all factors $1-s$ in the first parametrization and $1-t$ in the second. Since these factors are nonnegative on the square, nonnegativity of a quotient implies nonnegativity of the original polynomial.

	Rescale each dyadic subrectangle to $[0,1]^2$ with coordinates $(\xi,\eta)$. For a polynomial $\sum_{p=0}^m\sum_{q=0}^n a_{pq}\xi^p\eta^q$ of bidegree $(m,n)$, the power-to-Bernstein conversion~\cite[Sec.~3, eq.~(9)]{TitiGarloff2019} gives the following coefficient at indices $0\le i\le m$, $0\le j\le n$,
	\[
		\sum_{p=0}^i\sum_{q=0}^j
		a_{pq}\frac{\binom{i}{p}\binom{j}{q}}{\binom{m}{p}\binom{n}{q}} .
	\]
	The univariate Bernstein basis functions $\binom{m}{i}\xi^i(1-\xi)^{m-i}$ are nonnegative and sum to one by the binomial theorem. Their tensor products do likewise, so each polynomial is bounded below by its smallest Bernstein coefficient~\cite[eq.~(10)]{TitiGarloff2019}. For a given $\kappa$, check that all Bernstein coefficients of the reduced polynomials for $b,d,E$ are nonnegative. The same must hold for at least one of the reduced polynomials for $\chi$ and $\Delta$. Accept a rectangle if both conditions hold for one $\kappa\in\{0,\tfrac12,1\}$. Bisect the remaining rectangles at their midpoints, alternating the coordinates in the orders $(z,s)$ and $(t,z)$ for the two parametrizations. Exact rational arithmetic gives $62$ accepted rectangles of depth at most $11$.


\begin{thebibliography}{99}
	\bibitem{BoydVandenberghe2004}
	S.~Boyd and L.~Vandenberghe. \href{https://doi.org/10.1017/CBO9780511804441}{\emph{Convex Optimization}}.
	Cambridge University Press, 2004.

	\bibitem{DingSunYe2014}
	C.~Ding, D.~Sun, and J.~J.~Ye. \href{https://doi.org/10.1007/s10107-013-0735-z}{First order optimality conditions for mathematical programs with semidefinite cone complementarity constraints}.
	\emph{Mathematical Programming} \textbf{147} (2014), 539--579.

	\bibitem{GareyJohnsonStockmeyer1976}
	M.~R.~Garey, D.~S.~Johnson, and L.~Stockmeyer. \href{https://doi.org/10.1016/0304-3975(76)90059-1}{Some simplified NP-complete graph problems}.
	\emph{Theoretical Computer Science} \textbf{1} (1976), 237--267.

	\bibitem{Higham1988}
	N.~J.~Higham. \href{https://doi.org/10.1016/0024-3795(88)90223-6}{Computing a nearest symmetric positive semidefinite matrix}.
	\emph{Linear Algebra and its Applications} \textbf{103} (1988), 103--118.

	\bibitem{HornJohnson2013}
	R.~A.~Horn and C.~R.~Johnson. \href{https://doi.org/10.1017/CBO9781139020411}{\emph{Matrix Analysis}}, 2nd ed.
	Cambridge University Press, 2013.

	\bibitem{HouPengDing2026}
	W.~Hou, Y.~Peng, and C.~Ding.
	\emph{An Inexact Augmented Lagrangian Method for Semidefinite Complementarity Optimization}.
	\href{https://arxiv.org/abs/2611.14885v1}{arXiv:2611.14885v1}, 2026.

	\bibitem{HouPengDing2026Revised}
	W.~Hou, Y.~Peng, and C.~Ding.
	\emph{An Inexact Augmented Lagrangian Method for Semidefinite Complementarity Optimization with Applications}.
	arXiv:2611.14885, 2026.

	\bibitem{MooreKearfottCloud2009}
	R.~E.~Moore, R.~B.~Kearfott, and M.~J.~Cloud. \href{https://doi.org/10.1137/1.9780898717716}{\emph{Introduction to Interval Analysis}}.
	Society for Industrial and Applied Mathematics, 2009.

	\bibitem{TemmeDLMF}
	N.~M.~Temme. \href{https://dlmf.nist.gov/3}{Numerical methods}.
	Chapter~3 in \emph{NIST Digital Library of Mathematical Functions}, National Institute of Standards and Technology. Accessed 22 November 2026.

	\bibitem{TitiGarloff2019}
	J.~Titi and J.~Garloff. \href{https://doi.org/10.1016/j.amc.2018.08.049}{Matrix methods for the tensorial Bernstein form}.
	\emph{Applied Mathematics and Computation} \textbf{346} (2019), 254--271.

	\bibitem{WuZhangZhangLu2016}
	J.~Wu, Y.~Zhang, L.~Zhang, and Y.~Lu. \href{https://doi.org/10.1142/S0217595916500251}{A sequential convex program approach to an inverse linear semidefinite programming problem}.
	\emph{Asia-Pacific Journal of Operational Research} \textbf{33} (2016), 1650025.
	\end{thebibliography}
\end{document}